\documentclass[12pt]{article}
\usepackage{amsmath}
\usepackage{amssymb}

\newtheorem{theorem}{Theorem}

\def \bthm {\begin{theorem}}
\def \ethm {\end{theorem}}
\newtheorem{defi}{Definition}

\newtheorem{prop}{Proposition}
\newtheorem{rmk}{Remark}
\newtheorem{cor}{Corollary}

\def \bsk {\bigskip}
\def \msk {\medskip}

\def \qed {\hfill $\Box$}
\newcommand{\pff}[1]{\noindent{\bf Proof of #1. }}

\def \cir {\chi_1}
\def \omi {\omega_1}
\def \ali {\alpha_1}

\def \cF {\mathcal{F}}
\def \cG {\mathcal{G}}

\def \knk {KG(n,k)}

\begin{document}

\title{The robust chromatic number of certain graph classes}
\author{G\'abor Bacs\'o$^a$\and
  Csilla Bujt\'as$^{b,c}$\and
  Bal\'azs Patk\'os$^d$\and
  Zsolt Tuza$^{d,e}$  \qquad
  M\'at\'e Vizer$^f$ \\
  \small $^a$ Institute for
  Computer Science and Control \\
\small $^b$ Faculty of Mathematics and Physics, University of Ljubljana, Slovenia \\
\small $^c$ Institute of Mathematics, Physics, and Mechanics, Ljubljana, Slovenia \\
\small $^d$ Alfr\'ed R\'enyi Institute of Mathematics \\
\small $^e$ University of Pannonia \\
\small $^f$ Budapest University of Technology and Economics
  }
\date{
}
\maketitle

\begin{abstract}
    A 1-selection $f$ of a graph $G$ is a function $f: V(G)\rightarrow E(G)$ such that $f(v)$ is incident to $v$ for every vertex $v$. The 1-removed $G_f$ is the graph $(V(G),E(G)\setminus f[V(G)])$. The (1-)robust chromatic number $\chi_1(G)$ is the minimum of $\chi(G_f)$ over all 1-selections $f$ of $G$. 
    
    We determine the robust chromatic number of complete multipartite graphs and Kneser graphs and prove tight lower and upper bounds on the robust chromatic number of chordal graphs and some of their extensively studied subclasses, with respect to their ordinary chromatic number.
\end{abstract}

\section{Introduction}
Graph colorings and independent sets are central notions in graph theory. Various versions of graph colorings have been studied in the past decades. The focus of the present paper is the following recent variant.

\begin{defi}
For every nonnegative integer\/ $s$, an\/ \emph{$s$-selection}
 on\/ $G=(V,E)$ is an assignment\/ $f:V\to 2^{E}$
  such that\/ $f(v)\subseteq E(v)$ and\/ $|f(v)|\le s$, where\/ $E(v)$
 denotes the set of edges incident with\/ $v$.
 The graph\/ $G_f$ with vertex set\/ $V(G_f)=V(G)$ and edge set
 $$
   E(G_f) := E(G) \setminus \bigcup f(V(G)) 
 $$
  is termed an\/  \emph{$s$-removed subgraph} of\/ $G$.
Then\/ 
\begin{itemize}
    \item 
    the \textit{$s$-robust chromatic number of $G$ is} $\chi_s(G)=\min_f\chi(G_f)$,
    \item
    the \textit{$s$-robust independence number of $G$ is} $\alpha_s(G)=\max_f\alpha(G_f)$,
    \item
    the \textit{$s$-robust clique number of $G$ is} $\omega_s(G)=\min_f\omega(G_f)$,
\end{itemize}  
where $\min$ and $\max$ are taken over all $s$-selections of $G$.
\end{defi}

Observe that the ordinary chromatic, independence, and clique numbers of $G$ are $\chi_0(G)$, $\alpha_0(G)$, and $\omega_0(G)$, respectively. The notion of 1-robust chromatic number was introduced in \cite{PTV} as a tool to investigate specific Tur\'an-type problems. The systematic study of 1-robust parameters was initiated in \cite{BBPTV}. In this paper, we still concentrate on the $s=1$ case and say robust instead of 1-robust. Analogously, a coloring $c$ of $V(G)$ is robust if there exists a 1-selection $f$ such that $c$ is proper on $G_f$, and a subset $U\subset V(G)$ is robust independent if there exists a 1-selection $f$ such that $U$ is independent in $G_f$.


 



\subsection{New results}


A graph is \emph{chordal} if it does not contain any
 \emph{induced} cycles longer than~3.
Two famous subclasses (incomparable to each other) are the classes of 
 interval graphs and split graphs, the latter includes the
 proper subclass of threshold graphs (see definition in Section \ref{secchord}).
Here we establish tight inequalities for these graph classes.

The next theorem shows that the general lower
 bound $\cir(G) \geq \left\lceil \frac{\chi(G)}{3}\right\rceil$,
 which is valid for all graphs \cite{BBPTV}, is actually
  tight for an infinite subclass of threshold graphs, and a slightly weaker upper bound is valid for the more general class of \emph{split graphs}. We say that $G$ is \emph{$\omega$-unique} if it contains only one clique of order $\omega(G)$.

\msk

\begin{theorem} \label{t:threshold}
For every threshold graph\/ $G$,
$$
	\cir(G) = \left\{ 
	\begin{array}{ll} 
	\frac{\chi(G)}{3}+1,&\text{if}\ \chi(G) \equiv 0\ (\bmod\ 3) \mbox{ and\/ $G$ is not\/ $\omega$-unique,} \\ \left\lceil \frac{\chi(G)}{3}\right\rceil\!,& \mbox{otherwise.} 
	\end{array} \right.  
$$
For split graphs, the upper bound\/
 $\cir(G) \leq \left\lceil \frac{\chi(G)-1}{3}\right\rceil +1$
  is valid and tight, except for bipartite\/~$G$.
In particular, if\/ $G$ is a split graph with\/
 $\chi(G) \equiv 1\ (\bmod\ 3)$ then\/ $\cir(G)=\frac{\chi(G)+2}{3}$.
\end{theorem}

It is not true for the more general class of chordal graphs that $\chi(G)/3$ is an asymptotically tight upper bound.
Instead, the following holds.

\bthm\label{chordal}
$(i)$ \ If\/ $G$ is a chordal graph, then\/
   $$\omi(G) \leq \cir(G) \leq \left\lceil \frac{\chi(G)}{2} \right\rceil.$$

\noindent
$(ii)$ \ For every\/ $k\geq 2$ there exists an interval graph\/ $G_k$ such that\/ 
$\omega(G_k)=\chi(G_k)=k$ and\/
   $\omi(G_k) = \cir(G_k) = \left\lceil \frac{\chi(G_k)}{2} \right\rceil=\left\lceil \frac{k}{2} \right\rceil$.
\ethm

On the other hand, a further restriction on interval graphs drops $\cir$ down near $\chi/3$.

\bthm\label{unit}
There exists a constant $c$ such that for every unit interval graph $G$ we have $$\frac{\chi(G)}{3}\le \chi_1(G)\le \frac{\chi(G)}{3}+c.$$
\ethm

In Section 3 we solve the problem of determining $\cir$ for
 complete multipartite graphs $K_{n_1,\dots,n_t}$.
Throughout, we denote the number of vertex classes by $t$
 and write $n_i$ for the size of class $V_i$ for all $1\leq i\leq t$.
It will be assumed that the classes are in increasing order
 of their size, i.e.\ $n_1\leq n_2\leq \cdots \leq n_t$.
Then $\cir(K_{n_1,\dots,n_t})$ can be computed on the basis of the
 following two results, which complement each other.

\bthm   \label{t:cmk}
If\/ $n_t\leq 2$, assume that\/ $n_1=\ldots=n_p=1$ and
 $n_{p+1}=\ldots=n_{p+q}=2$, where\/ $p+q=t$.
Then
 $$
   \cir(K_{n_1,\dots,n_t}) =
    \left\lceil
     \frac{p + \left\lfloor 3q/2 \right\rfloor}{3} 
    \right\rceil .
 $$
\ethm

\bthm   \label{t:cmn}
If\/ $n_t\geq 3$, then an optimal 1-selection for\/ $\cir$ is
 obtained by enlarging\/ $V_t$ to an independent set with a
 vertex of\/ $V_1$.
That is,
 $$
   \cir(K_{n_1,\dots,n_t}) = 1 + \cir(K_{n_1-1,\dots,n_{t-1}}) .
 $$
In particular, if\/ $n_1=1$ and\/ $n_t\geq 3$, then
 $$
   \cir(K_{n_1,\dots,n_t}) = 1 + \cir(K_{n_2,\dots,n_{t-1}}) .
 $$
\ethm

In Section 4 we analyze the behavior of
 $\ali$ and $\cir$ in Kneser graphs.
Let $k\geq 2$ and $n\geq 2k$ be integers.
The \emph{Kneser graph} $\knk$ has $\binom{[n]}{k}=\{S\subseteq \{1,2\dots,n\}:|S|=k\}$ as its
 vertex set and two vertices are adjacent if and only if the
 corresponding $k$-element sets are disjoint.
Hence the \textit{intersecting subsystems}  of $\binom{[n]}{k}$ (those $\cF\subseteq \binom{[n]}{k}$ for which any $F,F'\in \cF$ have non-empty intersection) are in
 one-to-one correspondence with the independent sets of $\knk$.


\begin{theorem}\label{inters}
For any\/ $k\ge 2$ there exists\/ $n_0(k)$
 such that if\/ $n\ge n_0(k)$, then we have\/
 $$
 \alpha_1(KG(n,k))=\binom{n-1}{k-1}+1,
 $$ and\/
 $n_0(k)$ can be chosen to be\/ $8k^2$.
Furthermore, if $\cF$ is a\/ robust independent family
 in\/ $KG(n,k)$ such that for every \/ $x\in [n]$
 there exist at least two sets\/ $F_x,G_x\in \cF$
 with\/ $x\notin F_x,G_x$, then\/
  $|\cF|\le 8k\binom{n-2}{k-2}$ holds. 
\end{theorem}

\begin{theorem}\label{chi1kneser}
For any fixed\/ $k\geq 2$, we have\/ $$\chi_1(KG(n,k))=n-\Theta(n^{1/k})$$ as\/ $n\to\infty$.
\end{theorem}

 We finish the introduction by stating two propositions from \cite{BBPTV} that we will use in our proofs.
A graph is \textit{quasi-unicyclic} if each of its components contains at most one cycle.

\begin{prop}[\cite{BBPTV}]   \label{p:basic}
$(i)$ The value\/ $\cir(G)$ of a graph\/ $G=(V,E)$ is equal
 to the minimum number\/ $k$ of vertex classes in a partition\/
 $V=V_1\cup\cdots\cup V_k$ such that each\/ $V_i$ induces a
 quasi-unicyclic subgraph in\/ $G$.

$(ii)$ A graph\/ $G$ satisfies\/ $\cir(G)=1$
 if and only if it is quasi-unicyclic.
In particular, every tree has\/ $\cir=1$.

\end{prop}

\begin{prop}[\cite{BBPTV}] \label{f:chi}
For every graph\/ $G$ we have
\begin{equation} \label{eq:chi}
 \left\lceil \frac{\chi(G)}{3} \right\rceil \le \cir(G) \le \chi(G).
\end{equation}
All these bounds are tight, for all possible values of\/
 $\chi$.
\end{prop}

\section{Classes of chordal graphs}\label{secchord}

In this section we prove Theorem \ref{t:threshold}, Theorem \ref{chordal}, and Theorem \ref{unit}. First, we need to define the graph classes in Theorem \ref{t:threshold}.

A graph is called a \textit{split graph} if its vertex set can be partitioned into two sets, say $A$ and $B$, such that $A$ induces a complete subgraph and $B$ is independent.

A graph $G$ is called a \emph{threshold graph} if there exists a threshold $h$ and a function $g: V(G) \to \mathbb{R}$ such that $x, y\in V(G)$ are adjacent if and only if $g(x)+g(y) > h$. It follows from the definition (in fact, it is equivalent to the definition) that the vertex set of a threshold graph $G$ can be partitioned into two sets $A$ and $B$ (one of them might be empty) which satisfy the following conditions:
\begin{itemize}
\item[$(i)$] $A=\{a_1, \dots , a_q\}$ induces a maximum clique in $G$ and its vertices can be ordered such that $N[a_1] \supseteq N[a_2] \supseteq  \dots \supseteq N[a_q]$;
\item[$(ii)$] $B=\{b_1, \dots , b_s\}$ is an independent set in $G$ and its vertices can be ordered such that $N(b_1) \supseteq N(b_2) \supseteq  \dots \supseteq N(b_s)$;
\item[$(iii)$] $N(a_q)\cap B=\emptyset$.
\end{itemize}
A partition $(A, B)$ of $V(G)$ will be called a \emph{threshold partition} of $G$ if it satisfies $(i)-(iii)$. 

\msk

\pff{Theorem \ref{t:threshold}} We consider the threshold partition $(A,B)$ of $G$ with the notation introduced in $(i)-(iii)$. By definition, $A$ is a clique, $B$ is an independent set, and there is no clique that contains both $a_q$ and a vertex from $B$. It follows that $\omega(G)=q$ and, as every threshold graph is perfect, $\chi(G)=q$. Observe that $G$ is $\omega$-unique if and only if $a_{q-1}b_1$ is not an edge in $G$. Let $k=\lceil q/3\rceil -1$.

Suppose first that $q \not\equiv 0\ (\bmod\ 3)$, and define the vertex classes $V_j=\{a_{3j-2}, a_{3j-1}, a_{3j}\}$ for every $j\in [k]$.
The remaining vertices form the class $V_{k+1}= V(G) \setminus \bigcup_{j\in [k]}V_j$.
Observe that $G[V_j]$ is a unicyclic graph (in fact a $3$-cycle) for every $ j \in [k]$.
If $q \equiv 1\ (\bmod\ 3)$, then $V_{k+1}=\{a_q\} \cup B$ which is an independent set.
If $q \equiv 2\ (\bmod\ 3)$, then $V_{k+1}=\{a_{q-1}, a_q\} \cup B$ which induces a cycle-free graph as every edge of $G[V_{k+1}]$ is incident to $a_{q-1}$.
By Proposition~\ref{p:basic}$(i)$, the partition $V_1, \dots ,V_{k+1}$ defines a robust coloring for $G$ in both cases.
This implies $\cir(G) \leq k+1= \lceil \chi(G)/3 \rceil$ and, by the lower bound in (\ref{eq:chi}), we may conclude $\cir(G)= \lceil \chi(G)/3 \rceil$. 

Suppose now that $q \equiv 0\ (\mathrm{mod}\ 3)$ and $G$ is $\omega$-unique. For this case, we define $V_{k+1}=\{a_{q-2}, a_{q-1}, a_q\} \cup B$ and keep the notation $V_j=\{a_{3j-2}, a_{3j-1}, a_{3j}\}$ for $j\in [k]$.  Property $(iii)$ from the definition ensures $N[a_q]\cap B=\emptyset$. The $\omega$-uniqueness implies $a_{q-1}b_1 \notin E(G)$ that, together with property $(ii)$ gives $N[a_{q-1}] \cap B=\emptyset$. As $B$ is independent, every edge in $G[V_{k+1}]$ except $a_{q-1}a_q$ is incident to $a_{q-2}$ and therefore, $V_{k+1}$ induces a unicyclic graph. Since  $G[V_{j}]$ is also unicyclic for every $j \in [k]$, we infer that $V_1, \dots ,V_{k+1}$ gives a robust coloring for $G$ and $\cir(G) \leq \lceil \chi(G)/3 \rceil$. By inequality $(\ref{eq:chi})$, we conclude $\cir(G)= \lceil \chi(G)/3 \rceil$.

In the last case, $q \equiv 0\ (\mathrm{mod}\ 3)$ and $G$ is not $\omega$-unique. Observe first that $V_j=\{a_{3j-2}, a_{3j-1}, a_{3j}\}$ for $j\in [k+1]$ together with $V_{k+2}=B$ defines a robust coloring for $G$ and therefore, $\cir(G) \leq k+2$.  
Assume now for a contradiction that $W_1, \dots W_{k+1}$ is a robust color partition of $V(G)$. 
As $W_j$ induces a quasi-unicyclic graph and $A$ is a clique, $|W_j \cap A| \leq 3$ holds for every $j \in [k+1]$.
In fact, $|A|=3k+3$ implies $|W_j \cap A|=3$ for every color class. Since $G$ is not $\omega$-unique, $a_{q-1}b_1$ is an edge and $b_1$ is adjacent to all vertices but $a_q$ from $A$. Therefore, if $b_1$ is contained in the color class $W_i$, then $b_1$ is adjacent to at least two vertices from $W_i\cap A$ and $G[(W_i\cap A)\cup \{b_1\}]$ is either a complete graph $K_4$ or a $K_4-e$. In either case, $G[W_i]$ is not quasi-unicyclic, which is a contradiction. Thus, $\cir(G) >k+1$ holds and we conclude $\cir(G)=k+2= \frac{\chi(G)}{3}+1$. 

\msk

Split graphs also enjoy the property $\chi(G)=\omega(G)$.
Assume that $A\subset V(G)$ induces a complete subgraph of
 cardinality $\chi(G)$, and $B=V(G)\setminus A$ is an
 independent set.
Then we can have a robust coloring on $G[A]$ with
 $\left\lceil \frac{\chi(G)}{3}\right\rceil$ colors, and make
 $B$ monochromatic with a new color.
This $\left\lceil \frac{\chi(G)}{3}\right\rceil +1$ is the same
 as $\left\lceil \frac{\chi(G)-1}{3}\right\rceil +1$ unless
 $\chi(G)=3k+1$ for some integer $k$.
In that case we can have a robust $k$-coloring on $G[A]-v$ for
 a $v\in A$, and since $B\cup \{v\}$ induces a star, all edges
 from $B$ to $v$ can be omitted by a 1-selection.

The formula $\left\lceil \frac{\chi(G)-1}{3}\right\rceil +1$ is
 not tight if $G$ is a bipartite split graph, because in that case $G$ must be
 a double star, i.e.\ a particular tree, hence $\cir(G)=1$.
For larger $\chi=t\geq 3$, however, we can obtain a tight construction by taking $|A|=|B|=t$
 and putting $K_{t,t}-tK_2$ between $A$ and $B$.
For $t=3k+1$, we have seen that a robust coloring with
 $k+1$ colors is possible, and already the set $A$ requires
 that many colors.
For $t=3k$ and $t=3k-1$, the verified upper bound is $k+1$, and
 we argue that $k$ colors do not suffice (except if $3k-1=2$).

If $t=3k$, the only way for a robust $k$-coloring on $A$ is to
 select the edges of $k$ disjoint triangles in $G[A]$.
This cannot be extended to a robust $k$-coloring of $G$, because
 every $v\in B$ has at least two neighbors in each selected
 triangle in $A$, while only one incident edge
  can be deleted from $v$.

If $t=3k-1$, where $k\geq 2$, then a robust $k$-coloring on $A$
 has $k-1$ classes of size 3 (omitted triangles) and one class
 of size 2, say omitting the edge $xy$ by selecting $f(x)=xy$.
As in the previous case, the triangle classes do not admit any
 extension with vertices from $B$.
The 2-element class $\{x,y\}$ can be extended with the
 non-neighbor of $x$, with the non-neighbor of $y$, and with one
 further vertex $v\in B$ by defining $f(v)=xv$ and $f(y)=yv$.
But there are at least two more vertices in $B$, hence a
 further color will necessarily be used.
\qed
\bsk

\pff{Theorem \ref{chordal}}
$(i)$ \
Let $G=(V,E)$ be a chordal graph with $\chi(G)=\omega(G)=k\geq 2$.
Consider a proper vertex $k$-coloring of $G$, with color classes
 $V_1,\dots,V_k$.
Since $G$ is chordal, the union of any two color classes induces a forest in $G$.
Thus, there exists a 1-selection $f$ such that in the 1-removal
 $G_f$ all the sets $V_1\cup V_2$, $V_3\cup V_4$, $\dots,$
 $V_{k-1}\cup V_k$ are independent.
Hence the chromatic number has been decreased by
 $\lfloor k/2 \rfloor$, as needed.

\msk

$(ii)$\quad
The assertion is trivial for $k=2$ as shown by $K_2$, and for $k=3$ we can take e.g.\
 $K_4-e$, one edge deleted from the complete graph of order 4,
 which is an interval graph with $\chi=3$ and $\omi=2$ because
 $|E(K_4-e)|>4=|V(K_4-e)|$.

Let $G_2=R^{even}_1=K_2$, $G_3=R^{odd}_1= K_4-e$. For $k=2t>2$, we define $G_k=R^{even}_t$ and for $k=2t+1>3$ we define $G_k=R^{odd}_t$ as the graph obtained by taking
 three vertex-disjoint copies of $R^{even}_{t-1}$ ($R^{odd}_{t-1}$) together with two
 universal vertices.
Formally this means
 $$
   R^{even}_t = 3R^{even}_{t-1} \oplus K_2 \hskip 0.5truecm
   R^{odd}_t = 3R^{odd}_{t-1} \oplus K_2 \hskip 0.5truecm 
 $$
  (where $\oplus$ is the complete join operation).
As $R^{even}_1$ and $R^{odd}_1$ are interval graphs, and disjoint union does not change $\omega$, $\chi$ and being an interval graph, while adding a new universal vertex increases $\omega$ and $\chi$ by one, but keeps the property of being an interval graph, we obtain by induction that $R^{even}_t$ and $R^{odd}_t$ are interval graphs with $\omega(R^{even}_t)=k=2t$ and $\omega(R^{odd}_t)=k=2t+1$.
We need to prove that $\omi(R^{even}_t)\geq t$ and $\omi(R^{odd}_t)\geq t+1$ hold. As the proofs are almost identical, we only consider the case of $k=2t+1$ and omit the odd superscript.

We set $V_1=V(R_1)$ and $V_t=\{x_t,y_t\}$, the latter being the
 set of the two universal vertices in $R_t$.
Consider an arbitrary 1-selection $f$ in $R_t$.
This $f$ defines only two edges in $f(V_t)$, hence there is a copy
 of $R_{t-1}$ toward which no edge is selected for $x_t$ and $y_t$.
Let $V_{t-1}=\{x_{t-1},y_{t-1}\}$ denote the set of the two universal
 vertices in this copy of $R_{t-1}$.
Inside this $R_{t-1}$ subgraph, there is a copy of $R_{t-2}$ toward which no edge
 is selected for $x_{t-1}$ and $y_{t-1}$.
And so on, finally we obtain $t$ sets $V_1,V_2,\dots,V_t$ such
 that $V_i=\{x_i,y_i\}$ for all $2\leq i\leq t$ and the set
  $V_1\cup \cdots \cup V_t$ induces a subgraph, say $H$, in $R_t$ which is
   isomorphic to $K_{2t+2}-e$. Further, for every $i$ with $2 \leq i \leq t$, the selected edges $f(x_i)$ and $f(y_i)$
  are either contained in $\bigcup_{j=i}^t V_j$ or do not belong to $E(H)$.
The proof will be done if we show $\omega(H_f)=t+1$.
\medskip 

Consider the graph $F$ with 
$$V(F)=V_1\cup \cdots \cup V_t \quad \mbox{and} \quad
 E(F)=f(V_1\cup \cdots \cup V_t) \cap E(H).$$ 
Choose three distinct vertices $x_1,y_1,z_1\in V_1$ such that
 \begin{itemize}
   \item $x_1y_1 \notin E(F)$,
   \item $x_1 \notin f(z_1)$.
 \end{itemize}
Such $x_1,y_1,z_1$ exist because $f$ can select at most four of
 the five edges induced by $V_1$.

We let $F'$ be the induced subgraph of $F$ obtained by the
 deletion of the 
 vertex in $V_1\setminus \{x_1,y_1,z_1\}$.
The following procedure implies $\alpha(F') \geq t+1$.
 \begin{itemize}
   \item We put $V^0=V(F')$ and $I^0=\emptyset$. As long as the $V^j$ is not empty, select a vertex $v$ from $V^j$
    of degree 0 or 1 in $F[V^j]$, add $v$ to $I^j$ to obtain $I^{i+1}$ and delete $v$ from $V^j$ together with its neighbor
    if it has one to obtain $V^{j+1}$. It is easy to see that a vertex $v$ of degree 0 or 1 always exists in $V_m\cap V^j$ where $m$ is the minimum $i$ for which $V_i\cap V^j\neq \emptyset$.
 \end{itemize}
The selected vertices obviously form an independent set in $F'$,
 hence they induce a complete graph in $H_f$.
$F'$ has $2t+1$ vertices, and in each step, we delete at most two vertices, therefore at least
 $t+1$ vertices are selected at the end.
Vertex $x_1$ can be selected first, and there is a feasible choice
 for the next selection until the entire $F'$ is eliminated.
\qed

\bsk

We finish this section by proving Theorem \ref{unit}. A \textit{unit interval graph} is a graph of which the vertices $v_1,v_2,\dots,v_n$ are labelled with reals $r_1,r_2,\dots, r_n$ such that $v_i$ is joined to $v_j$ if and only $|r_i-r_j|<1$. The $p^{\mbox{\scriptsize{th}}}$ power $G^p$ of graph $G$ has the same vertex set as $G$, and two vertices are connected by an edge if and only if their distance in $G$ is at most $p$.

We shall apply the following result proved first in \cite{FH}.
Later developments and further references are reported in the Introduction of \cite{Sou}.

\bthm[Fine, Harrop \cite{FH}]\label{power}
A $n$-vertex graph $G$ is a unit interval graph if and only if there exist $n'\ge n$ and $p\ge 1$ such that $G$ is an induced subgraph of $P_{n'}^p$.
\ethm

As a matter of fact, the exponent $p$ can be chosen to be $\omega(G)-1$, which is the same as $\chi(G)-1$.

\bsk

\pff{Theorem \ref{unit}}
The lower bound follows from Proposition \ref{f:chi}.
For the upper bound assume that $G\subseteq H=P_{n'}^p$ where $p=\chi(G)-1$.
If $p=1$, then $G$ is a linear forest and $\cir(G)=1$.
Let $\chi(G)=p+1=3k-r\geq 3$ with $r\in\{0,1,2\}$, and assume without loss of generality that $n'=3t$, $P$ being the path $v_1v_2\dots v_{3t}$.
Then each triplet $S_i=\{v_{3i-2},v_{3i-1},v_{3i}\}$ ($1\leq i\leq t$) induces a $K_3$ in $H$, whose edges can be taken as a 1-selection $f$.
This decomposes $H_f$ into the 3-element independent sets $S_1,\dots,S_t$.
Moreover, if $|i-j|>k$, then there is no edge between $S_i$ and $S_j$.
Consequently the sets $\displaystyle \bigcup_{i\equiv r (\mathrm{mod}~k+1)} S_i$ are independent for each $r=0,1,\dots,k$, so that $\cir(G)\leq k+1$ holds.
Since $k\leq\frac{1}{3}(\chi(G)+2)$, the theorem follows.
\qed

\section{Complete multipartite graphs}   \label{s:compmul}

In this section, we prove Theorem \ref{t:cmk} and Theorem \ref{t:cmn}.
Before the proofs of these results let us mention some of their
 consequences. Recall that we assume $n_1 \leq \cdots \leq n_t$ when we use the notation $K_{n_1,\dots,n_t}$.

\begin{cor}   \label{c:compmul}
If\/ $n_1\geq t$, then\/ $\cir(K_{n_1,\dots,n_t}) = t$.
\end{cor}

\begin{cor}
If\/ $n_1<t$, let\/ $j$ denote the largest integer such that\/
 $n_1+\ldots+n_j\leq t-j$.
Then\/ $\cir(K_{n_1,\dots,n_t}) \leq t-j$.
\end{cor}

As shown by many examples (and already by the complete graph $K_n$),
 the upper bound $t-j$ is often far from being tight.

Before the proofs of Theorems \ref{t:cmk} and \ref{t:cmn}
 we observe that in a complete multipartite graph
 three types of independent sets (and their subsets) can be created 
 by the removal of a 1-selection:
  \begin{itemize}
   \item[(a)] three vertices $v_{i_1},v_{i_2},v_{i_3}$ from
    three distinct classes $V_{i_1},V_{i_2},V_{i_3}$,
     hence deleting the edges of a $C_3$;
   \item[(b)] four vertices $v'_{i_1},v''_{i_1},v'_{i_2},v''_{i_2}$
    from two distinct classes $V_{i_1},V_{i_2}$,
     hence deleting the edges of a $C_4$;
   \item[(c)] one vertex $v_{i_1}$ from vertex class $V_{i_1}$
    together with another class $V_{i_2}$,
     hence deleting the edges of a star.
  \end{itemize}

\bsk

\pff{Theorem \ref{t:cmk}}
Let us denote by $f(p,q)$ the value of $\cir$ under the assumptions
 of the theorem.
It is obvious that $f(p,q) = 
    \left\lceil
     \frac{p + \left\lfloor 3q/2 \right\rfloor}{3} 
    \right\rceil$
 is indeed valid if $p+q\leq 2$, and also if $p=3$ and $q=0$.
For the remaining cases, we apply induction on the number of vertices.
In the next cases, we check how the selection of an independent set $S$
 of type (a), (b), or (c) modifies $f(p,q)$ depending on the sizes
 of vertex classes met by $S$.
Assuming that $f$ correctly expresses the value of $\cir$ for all
 combinations of $p',q'$ with $p'+2q'<p+2q$, we obtain the
 following recursions:
 \begin{itemize}
  \item[(a.1)] $(1,1,1) \longrightarrow 1 + f(p-3,q)$;
  \item[(a.2)] $(1,1,2) \longrightarrow 1 + f(p-1,q-1)$;
  \item[(a.3)] $(1,2,2) \longrightarrow 1 + f(p+1,q-2)$;
  \item[(a.4)] $(2,2,2) \longrightarrow 1 + f(p+3,q-3)$;
  \item[(b.1)] $(2,2) \longrightarrow 1 + f(p,q-2)$;
  \item[(c.1)] $(1,1) \longrightarrow 1 + f(p-2,q)$;
  \item[(c.2)] $(1,2) \longrightarrow 1 + f(p-1,q-1)$;
  \item[(c.3)] $(2,2) \longrightarrow 1 + f(p+1,q-2)$.
 \end{itemize}
From these formulas the following ones are relevant:
 $$
   1 + f(p-3,q) \,, \quad 1 + f(p-1,q-1) \,, \quad
   1 + f(p,q-2) \,, \quad 1 + f(p+3,q-3) \,.
 $$
In this list (a.3), (c3), and (c.1) do not appear because they are
 superseded by (b.1) and (a1), respectively, which are their
 alternatives also structurally.
Note further that
 $$
   f(p-3,q) + 1 = f(p,q) = f(p,q-2) + 1
 $$
  and the reduction (c.3) or (a.1) can always be applied, hence
 $f(p,q)$ is a general upper bound on $\cir$.
But it is also a lower bound because, in the other two cases, we have
 $$
   f(p+3,q-3) + 1 \geq f(p-1,q-1) + 1 \geq f(p,q) \,.
 $$
This can be verified by comparing the numerators, namely
 $(p+3)+\lfloor 3(q-3)/2\rfloor = (p-1)+4+\lfloor 3(q-1)/2\rfloor-3
  > (p-1) + \lfloor 3(q-1)/2\rfloor$
  and
 $(p-1)+\lfloor 3(q-1)/2\rfloor+3 = p+\lfloor (3q+1)/2\rfloor
  \geq p + \lfloor 3q/2 \rfloor$.
\qed

\bsk

Before proving Theorem \ref{t:cmn}, let us state results from \cite{PTV} on bipartite and complete tripartite graphs.

\bthm   \label{t:PTV}
$(i)$ \cite[Proposition 2.6.]{PTV}
The complete tripartite graph\/ $K_{r,s,t}$ with\/ $1 \le r \le s \le t$ and\/ $t\ge 2$ satisfies\/  $\cir(K_{r,s,t}) = 2$ if and only if\/  $r \le 2$;
 otherwise\/ $\cir(K_{r,s,t}) = \chi(K_{r,s,t}) = 3$.
 


$(ii)$ \cite{PTV} A bipartite graph\/ $F$ has\/ $\cir(F)=2$ (i.e.,\/ $\cir(F)=\chi(F)$) if and only if it contains a component with more edges than vertices.
\ethm

\pff{Theorem \ref{t:cmn}}
The assertion is obvious for $t=2$, and its validity
 is easily derived from Theorem \ref{t:PTV} for $t=3$.
Assuming $t\geq 4$, let $m=\cir(K_{n_1,\dots,n_t})$ and consider
 a partition $(X_1,\dots,X_m)$ of the vertex set
 $V=V_1\cup\cdots\cup V_t$ into the minimum number of subsets $X_i$
  that all become
 independent after the removal of a suitably chosen 1-selection.
The statement of the theorem is $X_m=\{v_1\}\cup V_t$ for a $v_1\in V_1$.
If this is not the case, then we modify $(X_1,\dots,X_m)$ to another partition
 $(X_1',\dots,X_m')$ where $X_m'=\{v_1\}\cup V_t$ will hold.

There can be five types of $X_i$ in the partition:
 \begin{itemize}
  \item[(a)] $X_i \subset V_j$ for some $1\leq j\leq t$;
  \item[(b)] $|X_i\cap V_j|=|X_i\cap V_k|=|X_i\cap V_l|=1$
   for some $1\leq j<k<l\leq t$;
  \item[(c)] $|X_i\cap V_j|=|X_i\cap V_k|=2$
   for some $1\leq j<k\leq t$;
  \item[(d)] $|X_i\cap V_j|=1$ and $X_i\subset V_k$
   for some $1\leq j<k\leq t$;
  \item[(e)] $X_i\subset V_j$ and $|X_i\cap V_k|=1$
   for some $1\leq j<k\leq t$.
 \end{itemize}

There are several possible immediate simplifications in these types.
If $(X_i\setminus V_j) \neq V_k$ in (d), we can extend $X_i$ to
 contain the entire set $V_k$ and omit the vertices of
  $V_k\setminus X_i$ from the other sets $X_{i'}$ that meet $V_k$.
A similar step applies to (e), and also to (a) that yields then
 $X_i'=V_j$.
In fact, option (e) can be eliminated because (d) removes $V_k$
 while (e) removes $V_j$---plus one element for each---and we have
 $|V_j|\leq |V_k|$, hence the optimum with (d) is at least as good
 as the optimum with (e).
In the sequel, we analyze further ways of simplifying a partition.

\msk

(1)\quad
Assume first that $X_m = V_k$; it means type (a).
If $k\neq t$, we modify $X_m$ to $V_t$, and
 replace $|V_k|$ vertices in the sets $X_i$ meeting $V_t$ with the
 vertices of $V_k$ in a way that their sizes remain unchanged.
After that, a vertex $v_1\in V_1$ can be added to the modified $X_m$
 and the proof is done.

\msk

(2)\quad
Assume next that $X_m\cap V_j=v_j$ and $X_m\setminus\{v_j\} = V_k$
 (that is, type (d) occurs).
If $k\neq t$, we modify $X_m$ to $(X_m\setminus V_k)\cup V_t$, and
 replace $|V_k|$ vertices in the sets $X_i$ meeting $V_t$ with the
 vertices of $V_k$, as in case (1).
This finishes the proof if $|V_j|=1$, because in that case $V_j$ can
 play the role of $V_1$.
Hence suppose that $|V_j|\geq 2$ holds.

\msk

(2.1)\quad
If an $X_{i'}$ of type (b) or (d) exists that contains a single
 vertex $v_1$ from $V_1$, we switch the positions of $v_1$ and
 $v_j$; then $X_m$ is successfully modified to
  $X_m'=\{v_1\}\cup V_t$, and the proof is done.

\msk

(2.2)\quad
Otherwise, all $X_i$s meeting $V_1$ are of type (c).
Say, one of them is $X_i=\{v_1,v_1',v_i,v_i'\}$, where
 $v_1,v_1'\in V_1$ and $v_i,v_i'\in V_i$.
All the following subcases will lead to vertex partitions
 containing a class $\{v_1\}\cup V_t$ with $v_1\in V_1$.

\msk

(2.2.1)\quad
If there is a further $v_j'\in V_j$ and $X_{j'}=\{v_j'\}\cup V_k$
 of type (d), we replace $X_m$, $X_i$, and $X_{j'}$ with
 $\{v_1\}\cup V_t$, $\{v_i,v_i',v_j,v_j'\}$, and $\{v_1'\}\cup V_k$,
 respectively.

\msk

(2.2.2)\quad
If a $v_j'\in V_j$ is covered with $X_{j'}=\{v_j',v_k,v_l\}$ of type
 (b), we replace $X_m$, $X_i$, and $X_{j'}$ with
 $\{v_1\}\cup V_t$, $\{v_i,v_i',v_j,v_j'\}$, and $\{v_1',v_k,v_l\}$,
 respectively.

\msk

(2.2.3)\quad
If $V_j$ meets an $X_{j'}=\{v_j',v_j'',v_k,v_k'\}$ of type (c),
 we replace $X_m$, $X_i$, and $X_{j'}$ with
 $\{v_1\}\cup V_t$, $\{v_i,v_i',v_k,v_k'\}$, and $\{v_1'\}\cup V_j$,
 respectively.
This completes the proof in case (2).

\msk

From now on we can assume that the entire $V$ is partitioned into
 sets of types (b) and (c) only.

\msk

(3)\quad
If some $V_j$ meets more than one set of type (c), we can reduce
 the situation to case (1).
Indeed, say $X_i=\{v_j,v_j',v_k,v_k'\}$ and
 $X_{i'}=\{v_j'',v_j''',v_l,v_l'\}$.
Then we can replace these sets with $X_i'=V_j$ and
 $X'_{i'}=\{v_k,v_k',v_l,v_l'\}$.

\msk

(4)\quad
If $V_t$ meets a set $X_i=\{v_j,v_k,v_t\}$ of type (b) and a set
 $X_{i'}=\{v_l,v_l',v_t,v_t'\}$ of type (c), we can replace them
  with $X_i'=\{v_j\}\cup V_t$ of type (b) and
 $X_{i'}'=\{v_k\}\cup V_l$, hence reducing to case (2.2.1).

\msk

(5)\quad
If $V_t$ only meets sets of type (b), we take three of those sets
 say $X_{m-2},X_{m-1},X_m$.
Then we set $X_m' = V_t$, and split
 $(X_{m-2}\cup X_{m-1}\cup X_m) \setminus V_t$
  into two vertex triplets.
This reduces case (5) to case (1) and completes the proof of the theorem.
\qed

\begin{cor}
The graph invariant\/ $\cir$ is computable in polynomial time
 in the class of complete multipartite graphs.
\end{cor}

\section{Kneser graphs}   \label{s:knes}

In this section, we prove Theorem \ref{inters} and Theorem \ref{chi1kneser}.
Before presenting our results, we quote two fundamental theorems
 from extremal set theory that will serve as tools in our proofs.

\begin{theorem}[Erd\H os, Ko, Rado \cite{EKR}]\label{ekr}
For any integer\/ $k\geq 2$ and any\/ $n\ge 2k$,
 we have\/ $\alpha(KG(n,k))=\binom{n-1}{k-1}$.
\end{theorem}

\begin{theorem}[Hilton, Milner \cite{HM}]\label{hm}
For any\/ $n\ge 2k+1$, if\/
 $\cF\subseteq \binom{[n]}{k}$ is intersecting with\/
 $\bigcap_{F\in \cF}F=\emptyset$, then\/ $|\cF|\le \binom{n-1}{k-1}-\binom{n-k-1}{k-1}+1$.
\end{theorem}

\pff{Theorem \ref{inters}} 
Let $\cF$ be a 1-independent family of sets in $KG(n,k)$. We will use multiple times that $KG(n,k)[\cF]$ is $K_{2,3}$-free. If there exists a vertex $x$ such that there is at most one set $F\in \cF$ with $x\notin F$, then $|\cF|\le \binom{n-1}{k-1}+1$. So we can assume that for any $x$ there exist $F_x, G_x\in \cF$ with $x \notin F_x\cup G_x$. Then all but two sets containing $x$ must meet $F_x\cup G_x$, so $d_\cF(x)\le \binom{n-1}{k-1}-\binom{n-2k-1}{k-1}+2$ which is at most $ 2k\binom{n-2}{k-2}$ if $k\ge 3$.

Observe that $\cF$ must contain an intersecting family of size at least $\frac{1}{4}|\cF|$. Indeed, as $\cF$ is 1-independent, the degree sum in $KG(n,k)[\cF]$ is at most $2|\cF|$. We keep removing a maximum-degree set from $\cF$ as long as the maximum degree is larger than 1. In each step, the sum of the degrees decreases by at least 4, so in the end we have at least half of the sets in $\cF'$. Then $KG(n,k)[\cF']$ is a matching, so we can keep half of the sets to obtain an intersecting family. 

So let $\cF^*$ be a maximum-size intersecting subfamily of $\cF$. If $\bigcap_{F\in \cF^*}F=\emptyset$, then by Theorem \ref{hm}, we have $|\cF|\le 4|\cF^*|\le 4k\binom{n-2}{k-2}+4<\binom{n-1}{k-1}+1$ if $n\ge 4k^2$. If $\bigcap_{F\in \cF^*}F\neq \emptyset$, then $|\cF^*|$ is at most the maximum vertex-degree of $\cF$, which we showed above to be at most $2k\binom{n-2}{k-2}$ if $k\ge 3$. Then $|\cF|\le 4|\cF^*|\le 8k\binom{n-2}{k-2}\le \binom{n-1}{k-1}$, whenever $n\ge 8k^2$. Finally if $k=2$, then the result in the first paragraph states $d_\cF(x)\le \binom{n-1}{k-1}-\binom{n-2k-1}{k-1}+2=6$ and so $|\cF|\le 4|\cF^*|\le 24\le \binom{n-1}{k-1}$ if $n\ge 8\cdot 2^2\ge 25$.
\qed

\begin{rmk}
Observe that\/ $n_0(k)\ge 3k+1$ as\/ $\alpha_1(KG(3k,k))\ge \binom{3k-1}{k-1}+2$.
Indeed,\/ $\cF=\{F\in \binom{[3k]}{k}:1\in k\}\cup \{[k+1,2k],[2k+1,2k]\}$
 is\/ $1$-independent, as all sets in\/ $\{F\in \binom{[3k]}{k}:1\in k\}$
  but\/ $[k]$ intersect at least one o\/f $[k+1,2k],[2k+1,3k]$,
 and so\/ $KG(3k,k)[\cF]$ is a triangle with pendant edges from
 two vertices of the triangle.

Furthermore, if\/ $n=2k+m$ with\/ $1\le m<k$ and $f(m)$ denotes
 the maximum size of an intersecting family\/
  $\cG\subseteq \binom{[2k+m-1]}{k}$ with\/ $|G\cup G'|> k+m$ for all $G,G'\in \cG$, then\/
 $\alpha_1(KG(2k+m,k))\ge \binom{2k+m-1}{k-1}+f(m)$.
Indeed,\/ $\cF=\{F\in \binom{[2k+m]}{k}:2k+m\in F\}\cup \cG$
 is\/ $1$-independent if\/ $\cG$ is as above.
As\/ $\{F\in \binom{[2k+m]}{k}:2k+m\in F\}$ and\/ $\cG$ are
 intersecting,\/ $KG(2k+m,k)[\cF]$ is bipartite and for any\/
 $G,G'\in \cG$ there does not exist any\/ $k$-subset of\/
 $[2k+m]$ that is disjoint from both\/ $G,G'$ and so\/
  $KG(3k,k)[\cF]$ is a star forest.
\end{rmk}


\pff{Theorem \ref{chi1kneser}}
For the upper bound, observe that if for some integer $c$ we have $n-c \ge \binom{c}{k}-\binom{2k}{k}$, then $\chi_1(KG(n,k))\le n-c+1$ holds. Indeed, if we enumerate $\binom{[c]}{k}\setminus \binom{[2k]}{k}$ as $G_1,G_2,\dots, G_h$ with $h\le n-c$, then the families $\cG_i=\{G\in \binom{[n]}{k}: \max G=n+1-i\}\cup \{G_i\}$ are 1-independent as they induce a star in $KG(n,k)$ and they cover $\binom{[n]}{k}\setminus\binom{[2k]}{k}$. So adding $\binom{[2k]}{k}$ to the $\cG_i$s, we obtain a partition of $\binom{[n]}{k}$ into $n-c+1$ 1-independent families. As $k$ is fixed, we have $\binom{c}{k}-\binom{2k}{k}=\Theta(c^k)$, and so the largest value of $c$ for which  $n-c\ge \binom{c}{k}-\binom{2k}{k}$ holds is $\Theta(n^{1/k})$. This finishes the proof of the upper bound.

To prove the lower bound, we need the following definition. We say that a family $\cF$ of sets is \textit{star-like with center $x$} if $x$ is contained in all but at most one set $F$ of $\cF$. (So all $\cG_i$s in the previous paragraph are star-like.) Suppose that a robust coloring of $KG(n,k)$ contains $a$ star-like and $b$ non-star-like color classes. Observe that if $\cF$ is star-like with center $x$, then we can add any set $F$ containing $x$ to still have a star-like and therefore 1-independent family. Therefore, if $x_1,x_2,\dots,x_a$ are the centers of the star-like color classes, then we can assume that all non-star-like color classes $\cF$ are subfamilies of $\binom{[n]\setminus \{x_1,x_2,\dots,x_a\}}{k}$.

By Theorem \ref{inters}, all non-star-like color classes have 
size at most $8k\binom{n-a-2}{k-2}$. Therefore, we must have $$b\cdot 8k\binom{n-a-2}{k-2}+a\ge \binom{n-a}{k}.$$
If $a\ge n-2\cdot k!\cdot n^{1/k}$, then there is nothing to prove. Otherwise, $a\le \frac{1}{2}\binom{n-a}{k}$, and thus we must have 
\[
b \ge \frac{1}{16k}\frac{\binom{n-a-2}{k}}{\binom{n-a-2}{k-2}}\ge \frac{1}{32k^3}(n-a)^{2},
\]
if $n$ is large enough compared to $k$.
Therefore, the number $a+b$ of color classes is at least $a+\frac{1}{32k^3}(n-a)^2.$ This expression takes its minimum at $a=n-16k^3$ with value $ n-8k^3\ge n-n^{1/k}$.  
\qed

\section*{Acknowledgements}

This research was supported in part by the National Research, Development and Innovation Office --
NKFIH under the grants SNN 129364 and FK 132060.
Cs.\ Bujt\' as acknowledges the financial support from the Slovenian Research Agency under the core funding P1-0297.

\end{document}